\documentclass[12pt,a4paper]{article}

\usepackage[latin1]{inputenc}
\usepackage[english]{babel}
\usepackage{amsmath}\usepackage{latexsym}
\usepackage{textcomp}\usepackage{amsthm}
\usepackage{amssymb}
\usepackage{mathrsfs}

\newtheorem{teo}{Theorem}[section]
\newtheorem{coro}[teo]{Corollary}
\newtheorem{prop}[teo]{Proposition}
\newtheorem{lem}[teo]{Lemma}

\newtheoremstyle{normal}{\topsep}{\topsep}%
{}
{}
{\bfseries}
{.}
{ }
{\thmname{#1}\thmnumber{~#2}\thmnote{~(#3)}}
\theoremstyle{normal}

\newtheorem{remark}[teo]{Remark}

\newtheorem{algo}[teo]{Algorithm}

\author{Ludovic NAGESSEUR}
\title{A bundle method using two polyhedral approximations of the $\epsilon$-enlargement of a maximal monotone operator}

\renewcommand{\epsilon}{\varepsilon}
\renewcommand{\varnothing}{\emptyset}

\begin{document}
\maketitle
\renewcommand{\proofname}{\normalfont{\textbf{Proof}}}

\abstract{In this work, we develop a variant of a bundle method in order to find a zero of a maximal monotone operator. This algorithm relies on two polyhedral approximations of the $\epsilon$-enlargement of the considered operator, via a systematic use of the transportation formula. Moreover, the use of a double polyhedral approximation in our algorithm could inspire other bundle methods for the case where the given operator can be split as the sum of two other maximal monotone operators.}
\bigbreak


\section{Introduction and motivation}
Bundle methods have been widely employed in order to minimize nonsmooth convex functions (see for example \cite{MR2294039,MR1857057,MR1295240,MR2875584,MR2219129,MR2853160,MR1433963,MR1650317,MR1906232,MR2591628}). Those techniques are well known to be more implementable than the classical one; for example, they don't use resolvents of a maximal monotone operator as in  classical proximal point methods, or contrary to subgradients methods, the computation of a descent direction does not require the knowledge of the whole operator considered, and this, at each iteration. In the last 15 years, Burachik, Sagastiz\'abal, Svaiter and Solodov have developped in \cite{MR1737312,rbm}, two bundle methods involving a maximal monotone operator $T$ on $\mathbb{R}^N$ in order to solve inclusions of the following type:
\begin{equation}\label{inc1}
    0\in T(x).
\end{equation}
Those bundle methods are implementable versions of the algorithms studied in \cite{eemmop,hppa}.\bigbreak 

Given a maximal monotone operator $T$ on $\mathbb{R}^N$, we consider the problem of finding a solution of the inclusion (\ref{inc1}). A classical method to solve this problem is the Proximal Point Algorithm of Rockafellar (see \cite{rock}). It consists at each iteration to solve the inclusion: 
\begin{equation}\label{proxpb}
    0\in c_kT(x)+(x-x^k).
\end{equation}	 
An exact solution of (\ref{proxpb}) may be regarded as a pair $y^k,v^k$ satisfying
\begin{equation}\label{exppayv}
v^k\in T(y^k),\quad c_kv^k+y^k-x^k=0.
\end{equation}
The inexact proximal point algorithm of Rockafellar is based on the fact that  the next iterate $x^{k+1}$ is given by an approximate solution of (\ref{exppayv}), that is $x^{k+1}=y^k$, where 
\begin{equation}\label{eqapproxppa}
v^k\in T(y^k),\quad c_kv^k+(y^k-x^k)=e^k,
\end{equation}
with the error $e^k$ satisfying
\begin{equation}\label{conderrorppa}
\|e^k\|\leq\sigma_k, \quad \sum_{k=0}^{\infty}\sigma_k<\infty.
\end{equation}

There exist various other methods of proximal type to solve (\ref{inc1}). Most of them are implicit procedures, that is, they give only the existence of a sequence converging to a solution, but they are not implementable.\bigbreak

New "hybrid" variants of the proximal point algorithm have been recently  proposed. The basic idea of those methods is to combine at each iteration, an approximated solution of the  proximal problem (\ref{proxpb}), with a projection \cite{hppa}, or an extragradient step \cite{eemmop}. They are at the origin of bundle type algorithms for the resolution of the problem (\ref{inc1}), see \cite{rbm,MR1737312,MR2353163}. The $\epsilon$-\textit{enlargement} operator $T^{\epsilon}:\mathbb{R}^N\rightarrow 2^{\mathbb{R}^N}$ defined for all $x\in\mathbb{R}^N$ by
\begin{equation*}
   T^{\epsilon}(x):=\{x^*\in\mathbb{R}^N|\langle y^*-x^*,y-x\rangle\geq -\epsilon\,\,\,\forall (y,y^*)\in\mathrm{Gr}(T)\}, 
\end{equation*}
where $\mathrm{Gr}(T):=\{(z,z^*)\in\mathcal{H}\times\mathcal{H}|z^*\in T(z)\}$ is the graph of $T$, intervenes in all those bundle techniques, and is approached in some points by a polyhedral approximation.  

This operator has been first mentioned by J. Mart\'inez-Legaz and M. Théra in \cite{MR1422399}, but not studied. Then, a survey in finite dimension has been performed in \cite{emo}. Afterwards, in connection with the resolution of inclusion problems, an extension to Hilbert spaces has been proposed in \cite{MR1737312,eemmop}. Finally, it is in \cite{MR1716028} that the notion has been generalized  to Banach spaces.\bigbreak

Clearly, for a monotone operator $T:\mathbb{R}^N\rightarrow 2^{\mathbb{R}^N}$,
\begin{equation*}\label{t01}
    T^0(x)=\bigcap_{\epsilon>0}T^{\epsilon}(x),\quad\mbox{for all } x\in\mathbb{R}^N,
\end{equation*}
and $T$ is maximal monotone if and only if:
\begin{equation*}\label{te02}
    T^0(x)=T(x),\quad\mbox{for all } x\in\mathbb{R}^N.
\end{equation*}
It holds also that for any $x\in\mathcal{H}$ and $\epsilon\geq 0$, $T(x)\subset T^{\epsilon}(x)$; consequently, the $\epsilon$-enlargement of $T$ can be seen as an outer approximation of $T$.\bigbreak

The Hybrid Projection-Proximal Point Method (HPPM, \cite{hppa}) is based on (\ref{eqapproxppa}) with next iterate generated by
\begin{equation*}
 x^{k+1}=x^k-\displaystyle\frac{\langle x^k-y^k,v^k\rangle}{\|v^k\|^2}\,v^k,
\end{equation*}
and the error term of the method satisfies
\begin{equation}\label{errorhppm}
\|e^k\|\leq\sigma\max\{c_k\|v^k\|\,\mbox{;}\,\|y^k-x^k\|\},\quad\sigma\in[0,1), 
\end{equation}
Sagastiz\'abal and Solodov have proposed in \cite{rbm} a bundle method as a particular case of the HPPM. This algorithm preserves all the nice properties of convergence of the HPPM, but with the advantage that it is computationally implementable independently of the structure of $T$ and converges under no assumptions other than the maximal monotonicity of $T$ and $T^{-1}(0)\neq\varnothing$. In this bundle algorithm, the following choice is made:  $y^k=x^k-\sigma_ks^k$, where $\sigma_k>0$ and $s^k$ is obtained by a polyhedral approximation of $T^{\epsilon_k}(x^k)$ by using the transportation formula.
\bigbreak

Other proximal hybrid methods propose to take $v^k\in T^{\epsilon_k}(y^k)$ in (\ref{eqapproxppa}). We can cite for example the Hybrid Approximate Extragradient-Proximal Point Algorithm (HAEPPA, \cite{happa}), given by:
\begin{equation}\label{hya1+}
    v^k\in T^{\epsilon_k}(y^k),\quad 0=c_kv^k+(y^k-x^k)-e^k,
\end{equation}
\begin{equation}\label{hya2+}
    x^{k+1}=x^k-c_kv^k,
\end{equation}
\begin{equation}\label{hya3++}
    \|e^k\|^2\leq\sigma^2\|y^k-x^k\|^2-2c_k\epsilon_k,
\end{equation}
where $\sigma\in[0,1)$, $c_k>0$ for all $k$, $e_k$ being the error term of the method. The positive sequence $\{\epsilon_k\}$ is convergent to zero; thus, $T^{\epsilon_k}$ approaches $T$ when $k\rightarrow\infty$. It follows that the relation (\ref{hya1+}) is a perturbed version of (\ref{eqapproxppa}), and   $\{\epsilon_k\}$ can be then considered as another error term.  Observe that if the $k$th proximal point subproblem is solved exactly, then the $k$th iteration of the above algorithm coincides with the classical proximal point iteration (this is because $c_kv^k+(y^k-x^k)=0$ and $\epsilon_k=0$ imply that $x^{k+1}=y^k$ and $v^k\in T(y^k)$). Therefore, in the special case $\sigma=0$, we retrieve the classical (exact) proximal point method.

In \cite{MR1871872}, the authors adopt the same idea as for the HAEPPA, that is their algorithm satisfies (\ref{hya1+}); also, the next iterate is given by an iteration similar to (\ref{hya2+}), and the error verifies
\begin{equation}\label{errorppm}
    \|e^k\|^2\leq\sigma_k^2(\|c_kv^k\|^2+\|y^k-x^k\|^2)-2c_k\epsilon_k,
\end{equation}
where $\sigma_k<1$ for all $k$.
\bigbreak

The Implementable Bundle Strategy (IBS) developped in \cite{MR1737312} can be outlined as follows. Given an arbitrary $y$ and $v\in T(y)$, the monotonicity of $T$ implies that $T^{-1}(0)$ is contained in the halfspace 
\begin{equation}\label{defhalfspace}
H_{y,v}:=\{z\in\mathbb{R}^N:\langle z-y,v\rangle\leq 0\}.
\end{equation}
Then, given a current iterate $x^k\notin T^{-1}(0)$,
\begin{description}
    \item[-] first, find $y^k$ and $v^k\in T(y^k)$ such that $x^k\notin H_{y_k,v_k}$,
    \item[-] then, project $x^k$ onto $H_{y_k,v_k}\supseteq T^{-1}(0)$ to obtain a new iterate 
\begin{equation*}
 x^{k+1}=P_{H_{y^k,v^k}}(x^k)=x^k-\displaystyle\frac{\langle x^k-y^k,v^k\rangle}{\|v^k\|^2}\,v^k.
\end{equation*}
\end{description}

By an elementary property of orthogonal projections, $x^{k+1}$ is closer to $T^{-1}(0)$ than $x^k$. However, in order to have a significant progress from $x^k$ to $x^{k+1}$, adequate choices are made in the IBS: the vector $v^k\in T(y^k)$ being given by an oracle, only the selection of the vector $y^k$ is possible; in \cite{MR1737312}, they choose 
\begin{equation}\label{defykibs}
y^k=x^k-\sigma_ks^k,
\end{equation}
where $s^k=P_{T^{\epsilon_k}(x^k)}(0)$ and $\{\sigma_k\}$ is a sequence of positive reals. The vector $s^k$ is obtained by using the transportation formula and is then the solution of the convex problem of finding the element of minimal norm of a polyhedral approximation of $T^{\epsilon_k}(x^k)$.\bigbreak

In this paper, we propose a new bundle algorithm of proximal type, based on algorithm IBS, but using a double polyhedral approximation via a systematic use of the Transportation Formula: $\{y^k\}$ is constructed as in (\ref{defykibs}), and thus, a first polyhedral approximation of the $\epsilon$-enlargement of $T$ is used in the calculus of $s^k$; the second approximation intervenes when searching a good candidate for the vector $v^k$. This is the scheme of the method:
\begin{equation}\label{vk2bm}
    v^k \mbox{ is an approximate vector of } T^{\epsilon_k}(y^k),
\end{equation}
\begin{equation}\label{iterprox}
    c_kv^k+(y^k-x^k)-e^k=0,
\end{equation}
\begin{equation}\label{xikinct}
    \xi^k\in T(y^k),
\end{equation}
\begin{equation}\label{iteratexk2bm}
    x^{k+1}=P_{H_{y^k,\xi^k}}(x^k)=x^k-\frac{\langle x^k-y^k,\xi^k\rangle}{\|\xi^k\|^2}\xi^k,
\end{equation}
\begin{equation}\label{ineqarror2bm}
\|e^k\|^2\leq c_k^2\|v^k\|^2+\|y^k-x^k\|^2,
\end{equation}
 
The vector $\xi^k$ plays a double role in this algorithm: first, since it satisfies the inclusion (\ref{xikinct}), it enables to construct a bundle of information which will permit to obtain the vector $v^k$ in (\ref{vk2bm}) by the transportation formula, and $\xi^k$ intervenes also in the construction of the new iterate of the sequence $\{x^k\}$ in (\ref{iteratexk2bm}). The sequence $\{c_k\}$ is supposed to be positive. Thanks to the choice of the parameters in our algorithm, the error term satisfies (\ref{ineqarror2bm}). We can note that the condition (\ref{ineqarror2bm}) is weaker than either (\ref{errorhppm}), (\ref{hya3++}) or (\ref{errorppm}), because the right hand-side in (\ref{ineqarror2bm}) is larger: computationally, at each iteration, the condition (\ref{ineqarror2bm}) is easier to verify.
\bigbreak
This work was motivated by the creation of a bundle method for the HAEPPA. Indeed, as for our method, we could use a polyhedral approximation to obtain the descent direction $s^k$ and another to get the vector $v^k$. It is interesting to study a "bundle adaptation" for the HAEPPA because a modified forward-backward splitting method of Tseng \cite{mfb} can be seen as a particular case of this method (see \cite[Section 5]{happa}). We could then obtain the (first) bundle method coming from a splitting one. But the main difficulties are the choices of the vector $y^k$ and of the sequence $\{c_k\}$ (this latter must be bounded away from zero), and also how to establish the relation of the error (\ref{hya3++}).

\bigbreak

\section{Preliminary results}

\subsection{Notations and assumptions} 

Given a multifunction $F:\mathbb{R}^N\rightarrow 2^{\mathbb{R}^N}$ and a set $E\subseteq\mathbb{R}^N$:

\begin{description}
    \item[-] the \textit{closure} of $E$ is denoted by $\overline{E}$,
    \item[-] the \textit{domain} of $F$ is $\mathrm{D}(F):=\{x\in\mathbb{R}^N:F(x)\neq\varnothing\}$ and
    \item[-] the \textit{graph} of $F$ is $\mathrm{Gr}(F):=\{(x,y)\in\mathbb{R}^N\times\mathbb{R}^N|y\in F(x)\}$.
    \item[-] We define the set $F(E):=\bigcup_{e\in E}F(e)$.
    \item[-] $B(x,\rho)$ denotes the opened ball centered in $x$ with radius $\rho$.
    \item[-] $F$ is \textit{locally bounded} at $x$ if there exists a neighbourhood $U$ of $x$ such that the set $F(U)$ is bounded.
    \item[-] $F$ is monotone if $\langle u-v,x-y\rangle\geq 0$ for all $x,y\in\mathbb{R}^N$, and all $u\in F(x)$, $v\in F(y)$.
    \item[-] $F$ is \textit{maximal monotone} if it is monotone, and, additionally, whenever there is some monotone operator $T$ such that $F(x)\subset T(x)$ for all $x\in\mathbb{R}^N$, this implies $F=T$.
    \item[-] $F$ is \textit{firmly nonexpansive} if for all $x,y\in\mathbb{R}^N$,
\begin{equation*}
    \|F(x)-F(y)\|^2+\|(\mathrm{Id}-F)(x)-(\mathrm{Id}-F)(y)\|^2\leq\|x-y\|^2.
    \end{equation*}
\end{description}\bigbreak

Recall that every maximal monotone operator is locally bounded in the interior of its domain (\cite[Theorem 1]{MR0253014}).
All along this paper, we suppose that $T$ in (\ref{inc1})  is defined in the whole of $\mathbb{R}^N$, so that it maps bounded sets in bounded sets. Finally, recall that the solution set of (\ref{inc1}) is assumed to be not empty.\bigbreak

\subsection{The $\epsilon$-enlargement. Some useful properties}

\subsubsection{Continuity properties} 

According to \cite[Corollary 32.10]{MR1716028}, the $\epsilon$-enlargement is locally bounded on the interior of its domain. Also, in \cite[Theorem 2.6]{eemmop}, it is shown that $T^{\epsilon}$ is Lipschitz-continuous whenever $\epsilon>0$. More precisely:

\begin{teo}\label{teolipteps}
Let $T:\mathbb{R}^N\rightarrow 2^{\mathbb{R}^N}$ be a maximal monotone operator such that $\mathrm{int(D}(T))$ is nonempty. Let $K\subset \mathrm{int(D}(T))$ be a compact set, and $0<\underline{\epsilon}<\overline{\epsilon}<+\infty$. Then, there exist nonnegative constants $A$ and $B$ such that for any $(\epsilon_1,x^1), (\epsilon_2,x^2)\in[\underline{\epsilon},\overline{\epsilon}]\times K$ and $v^1\in T^{\epsilon_1}(x^1)$, there exists $v^2\in T^{\epsilon_2}(x^2)$ satisfying:
\begin{equation}\label{lipcontteps}
   \|v^1-v^2\|\leq A\|x^1-x^2\|+B|\epsilon_1-\epsilon_2|.
\end{equation}
\end{teo}\bigbreak

 We can add this following result which establish the closedness of the graph of $T^{\epsilon}$:

\begin{prop}\,\rm\cite[Proposition 1 (iv)]{emo}\it\label{clteps}
For any maximal monotone operator $T:\mathbb{R}^N\rightarrow 2^{\mathbb{R}^N}$ and any sequence $\{(\epsilon_i,x^i,w^i\in T^{\epsilon_i}(x^i))\}_i$ such that $\epsilon_i\geq 0$ for all $i$,
\begin{equation*}
    \lim_{i\rightarrow\infty}x^i=x,\quad \lim_{i\rightarrow\infty}\epsilon^i=\epsilon,\quad \lim_{i\rightarrow\infty}w^i=v\quad \Rightarrow\quad v\in T^{\epsilon}(x).
\end{equation*}
\end{prop}\bigbreak

\subsubsection{The Transportation Formula}\label{sectionft}

In this section, we will use the notation
\begin{equation*}
   \Delta_m:=\left\{\alpha=(\alpha_1,...,\alpha_m)\in\mathbb{R}^m|\forall i=1,...,m,\,\,\alpha_i\geq 0,\,\, \sum_{i=1}^m\alpha_i=1\right\}
\end{equation*}
to designate the unit simplex of $\mathbb{R}^m$.\bigbreak

\begin{teo}\rm\textbf{\cite[Theorem 2.3]{eemmop}}\it\label{teoft}\,\,
Let $T$ be a maximal monotone operator on $\mathbb{R}^N$. Consider the set of  $m$ triplets:
\begin{equation*}
    \{(\epsilon_i\geq 0,z^i\in\mathbb{R}^N,w^i\in T^{\epsilon_i}(z^i))\}_{i=1,...,m}.
\end{equation*}
For all $\alpha\in\Delta_m$, set:\\
\begin{center}
$\displaystyle\hat{x}:=\sum_{i=1}^{m}\alpha_iz^i$;\\
$\displaystyle\hat{u}:=\sum_{i=1}^{m}\alpha_iw^i$;\\
$\displaystyle\hat{\epsilon}:=\sum_{i=1}^{m}\alpha_i\epsilon_i+\sum_{i=1}^{m}\alpha_i\langle w^i-\hat{u},z^i-\hat{x}\rangle$.
\end{center}
Then, $\hat{\epsilon}\geq 0$, and $\hat{u}\in T^{\hat{\epsilon}}(\hat{x})$.
\end{teo}\bigbreak

\begin{remark}
Observe that when $\epsilon_i=0$, $\forall i=1,...,m$, this theorem shows how to construct $\hat{v}\in T^{\hat{\epsilon}}(\hat{x})$, by using convex  combinations of pairs $(z^i,w^i)\in \mathrm{Gr}(T)$. In other words, the set $\mathrm{conv}(\{(z^i,w^i\in T(z^i))\}_{i=1,...,m})$ constitutes a polyhedral approximation of $\mathrm{Gr}(T^{\hat{\epsilon}})$.
\end{remark}\bigbreak

The following result is a consequence of Theorem \ref{teoft} in the particular case where $\epsilon_i=0$ for all $i=1,...,m$. It gives us a bound for $\hat{\epsilon}$. It will be very useful in the sequel. We can find it in \cite[Corollary 2.3]{MR1737312} or in \cite[Corollary 5.5.7]{MR2353163}.\bigbreak

\begin{coro}\label{corotf}
Consider the same notations as in the above theorem. Suppose that $\epsilon_i=0$, $\forall i\leq m$. Let $\tilde{x}\in\mathbb{R}^N$ and $\rho>0$ such that:
\begin{equation*}
    \|z^i-\tilde{x}\|\leq\rho,\quad \forall i\leq m.
\end{equation*}
Then, the convex sum:
\begin{equation*}
    (\hat{x},\hat{u}):=\left(\sum_{i=1}^m\alpha_iz^i,\sum_{i=1}^m\alpha_iw^i\right)
\end{equation*}
satisfies:
\begin{center}
    $\displaystyle\|\hat{x}-\tilde{x}\|\leq\rho$,\\
    $\displaystyle\hat{u}\in T^{\hat{\epsilon}}(\hat{x})$, with $\displaystyle\hat{\epsilon}:=\sum_{i=1}^m\alpha_i\langle w^i-\hat{u},z^i-\hat{x}\rangle\leq 2\rho M$,
\end{center}
where $M:=\max\{\|w^i\||i=1,...,m\}$.
\end{coro}\bigbreak

\begin{remark}
Suppose that $\hat{\epsilon}>0$. By rewritting the inequality (\ref{lipcontteps}) with $x^1=\hat{x}$, $x^2=\tilde{x}$, $v^1=\hat{u}$ and $\epsilon_1=\epsilon_2=\hat{\epsilon}$, the last Corollary and the Lipschitz continuity of $T^{\hat{\epsilon}}$ give the existence of a constant $A>0$ and an element $\tilde{u}\in T^{\hat{\epsilon}}(\tilde{x})$ such that:
\begin{equation*}
   \|\hat{u}-\tilde{u}\|\leq A\|\hat{x}-\tilde{x}\|\leq A\rho.
\end{equation*}
Thus, if $\rho$ is small (i.e. converges to $0$), more precisely, if we choose the vectors $z^i$ ($i=1,...,m$) in a small neighborhood of $\tilde{x}$, we will have: $\tilde{u}\approx\hat{u}$. In other words, $\hat{u}\in T^{\hat{\epsilon}}(\hat{x})$ approaches a vector $\tilde{u}\in T^{\hat{\epsilon}}(\tilde{x})$.
\end{remark}\bigbreak

\section{A bundle method with two polyhedral approximations of the $\epsilon$-enlargement of a maximal monotone operator}

Our algorithm respects the scheme (\ref{vk2bm})-(\ref{ineqarror2bm}). Consequently, we will make a good choice for $y^k$, $v^k$ in order to verify all the relations. Like for the algorithms proposed in \cite{eemmop}, \cite{MR1737312} and \cite{rbm}, we will take here $y^k=x^k-\sigma_ks^k$, where $s^k$ is given by a projection of $0$ onto a polyhedral approximation of $T^{\hat{\epsilon}_k}(x^k)$, where $\{\hat{\epsilon}_k\}$ is a  sequence of positive reals. In order to obtain $s^k$, the algorithm use the transportation formula to constitute a set $Q_p$ such that the convex hull of a part of its elements gives a polyhedral approximation of $T^{\hat{\epsilon}_k}(x^k)$ at step 1. The choice of $y^k$ appears useful in the construction of $Q_p$. The set $Q_p$ is also used to approach a vector $v^k\in T^{\epsilon_k}(y^k)$.\bigbreak

We will note
\begin{equation*}
    \Delta_I:=\left\{\lambda=\{\lambda_i\}_{i\in I}:\forall i\in I,\, \lambda_i\in\mathbb{R}_+,\,\,\sum_{i\in I}\lambda_i=1\right\}
\end{equation*}
the unit simplex associated with the set of index $I$.

\bigbreak

\noindent\hrulefill

\begin{algo}\label{bmhaeppa}
\textbf{A bundle type algorithm with a double polyhedral approximation of the $\epsilon$-enlargement of a maximal monotone operator}\vspace{0.5cm}\\
\textbf{Data:} $x^0\in\mathbb{R}^N$, $\tau>0$, $R>0$.\vspace{0.5cm}\\
\textbf{Initialization:} $k:=0$, $p:=0$.\vspace{0.5cm}\\
\textbf{Step 0: Stopping Test}\\ 
\begin{description}
\item{(0.a)} Compute $u^k\in T(x^k)$; If $u^k=0$, then STOP.\\
\item{(0.b)} Else, $p:=p+1$, set $(z^p,w^p):=(x^k,u^k)$. Set $n:=0$.\\  
\end{description}
\textbf{Step 1: Computing Search Direction}\\
\begin{description}
\item{(1.a)} Set $j:=0$.\\
\item{(1.b)} Define $\widehat{I}_{k,n,j}:=\{1\leq i\leq p|\|z^i-x^k\|\leq R\,2^{-j}\}$.\\
\item{(1.c)} Compute $\alpha^{k,n,j}:=\mathrm{argmin}\left\{\left\|\displaystyle\sum_{i\in \widehat{I}_{k,n,j}}\alpha_iw^i\right\|^2|\alpha\in\Delta_{\widehat{I}_{k,n,j}}\right\}$.\\
\item{(1.d)} Take $s^{k,n,j}:=\displaystyle\sum_{i\in \widehat{I}_{k,n,j}}\alpha_i^{k,n,j}w^i$.\\
\item{(1.e)} If $\|s^{k,n,j}\|\leq \tau\, 2^{-j}$,
then set $j:=j+1$, and LOOP to (1.b).\\
\item{(1.f)} Else, define $j_{k,n}:=j$ and $s^{k,n}:=s^{k,n,j_{k,n}}$. 
\bigbreak
\end{description}\vspace{0.5cm}
\textbf{Step 2: Line Search}\\
\begin{description}
\item{(2.a)} Set $l:=0$.\\
\item{(2.b)} Define 
$\sigma_{k,n,l}:=R\, 2^{-l}/\|s^{k,n}\|$.\\ 
\item{(2.c)} Define $y^{k,n,l}:=x^k-\sigma_{k,n,l}\,s^{k,n}$ and take $\xi^{k,n,l}\in T(y^{k,n,l})$.\\
\item{(2.d)} Define $I_{k,n,l}:=\{1\leq i\leq p|\|z^i-y^{k,n,l}\|\leq R\,2^{-l}\}$.\\ 
\item{(2.e)} Take $\lambda\in\Delta_{I_{k,n,l}}$ and set $v^{k,n,l}:=\displaystyle\sum_{i\in I_{k,n,l}}\lambda_iw^i$.\\

\item{(2.f)} 
Set
\begin{equation*}
    \epsilon_{k,n,l}:=\displaystyle\sum_{i\in I_{k,n,l}}\alpha^{k,n,j_{k,n}}_i\bigg\langle w^i-v^{k,n,l},z^i-\sum_{m\in I_{k,n,l}}\lambda_mz^m\bigg\rangle.
\end{equation*}
\item{(2.g)} If\begin{equation}\label{vsstrict}
    \langle v^{k,n,l},s^{k,n}\rangle < \frac{1}{2}\|s^{k,n}\|^2
\end{equation}
\hspace{6.5cm} or\\
    \begin{equation}\label{vxistrict}
    \langle s^{k,n},\xi^{k,n,l}\rangle<\frac{1}{2}\,\|s^{k,n}\|^2, 
\end{equation}
 \\ and ($l<j_{k,n}+1$), then\\ Set $l:=l+1$ and LOOP to (2.b).\\
\item{(2.h)} Else, define $l_{k,n}:=l$ and $y^{k,n}:=y^{k,n,l_{k,n}}$, $v^{k,n}:=v^{k,n,l_{k,n}}$, $\xi^{k,n}=\xi^{k,n,l_{k,n}}$, $\sigma_{k,n}=\sigma_{k,n,l_{k,n}}$, $\epsilon_{k,n}=\epsilon_{k,n,l_{k,n}}$.
\end{description}\vspace{0.5cm}
\textbf{Step 3:}\\
\begin{description}
\item{(3.a)} If \hspace{6cm}\textit{Null Step}\\ \begin{equation*}
    \langle v^{k,n},s^{k,n}\rangle < \frac{1}{2}\|s^{k,n}\|^2
\end{equation*}
 \hspace{6.5cm} or\\
    \begin{equation*}
    \langle s^{k,n},\xi^{k,n}\rangle<\frac{1}{2}\|s^{k,n}\|^2, 
\end{equation*}then\\
\\ Set $p:=p+1$, $(z^p,w^p):=(y^{k,n},\xi^{k,n})$. Set $n:=n+1$ and LOOP to (1.b).\\
\item{(3.b)} Else,\hspace{6cm} \textit{Serious Step}\\ Define $n_k:=n$, $j_k:=j_{k,n_k}$, $s^k=s^{k,n}$, $y^{k}=y^{k,n}$, $v^k:=v^{k,n_k}$, $\xi^k:=\xi^{k,n_k}$, $\sigma_k=\sigma_{k,n}$, $\epsilon_k:=\epsilon_{k,n}$.\\
Compute $$x^{k+1}:=x^k-\frac{\langle x^k-y^k,\xi^k\rangle}{\|\xi^k\|^2}\xi^k.$$\\
Set $k:=k+1$ and LOOP to Step 0.
\end{description}
\end{algo}

\noindent\hrulefill

\bigbreak

In the algorithm, we execute a double task simultaneously: the generation of the sequence $\{x^k\}$ and the construction of the bundle $\Gamma_p:=\{(z^0,w^0),...,(z^p,w^p)\}$. This double task comes from the \textit{serious} and the \textit{null} steps. The first, indexed in $k$, produces a new iterate $x^{k+1}$, while in the second, indexed in $n$, the pair $(y^{k,n},\xi^{k,n})\in\mathrm{Gr}(T)$ is renamed $(z^p,w^p)$ and is then added to the bundle.\bigbreak

At each step, the convex hull of  $Q_p:=\{w^0,..., w^p\}$ enables to obtain a polyhedral approximation of $T^{\hat{\epsilon}_k}(x^k)$ and of $T^{\epsilon_k}(y^k)$. It is first used  to generate a direction $s^k$ at step 1. During the steps 0 and 3, $Q_p$ grows; then, in order to limit the calculus and have a better approximation of $T^{\hat{\epsilon}_k}(x^k)$ and  $T^{\epsilon_k}(y^k)$, we use only the $w^i$ associated with vectors $z^i$ which are in a neighbourhood of $x^k$ or of $y^k$. Then, we can distinguish three bundles: the \textit{raw bundle} $\Gamma_p$, and two \textit{reduced bundles} $\{(z^i,w^i):i\in \widehat{I}_{k,n,j}\}$ and $\{(z^i,w^i):i\in I_{k,n,l}\}$. Actually, it is those two reduced (sub-) bundles which are used by the method to approach respectively $T^{\hat{\epsilon}_{k,n,j}}(x^k)$ at step 1, and $T^{\epsilon_{k,n,l}}(y^k)$ at each loop in $l$ at step 2 (for more precisions, see Remark \ref{remapp}).\bigbreak

Let us describe now the different steps of Algorithm \ref{bmhaeppa}:\bigbreak

\begin{itemize}
    \item The step 0 is the stopping test; an oracle gives an element $u^k\in T(x^k)$, and if $u^k=0$, then $x^k$ is a solution of the inclusion (\ref{inc1}).
 
    \item At step 1, we search a direction $s^k$. The set $T^{\hat{\epsilon}_{k,n,j}}(x^k)$ is replaced here by its polyhedral approximation formed by the convex hull of some selected elements in the range of $T$. The construction of $Q_p$ is done only with vectors $w^i\in T(z^i)$ such that $z^i$ belongs to a ball $\overline{B(x^k,R\,2^{-j})}$. The radius of this ball is successively divided by $2$ through a loop in $j$ in this step 1, until to the direction of minimal norm $s^{k,n,j}$ at step (1.d) has a norm bigger than $\tau 2^{-j}$. As the radius of the ball is reduced at each loop, the size of the set $\widehat{I}_{k,n,j}$ (that is the cardinal of the vectors defining the bundle $\{(z^i,w^i):i\in \widehat{I}_{k,n,j}\}$) is reduced.

\item At step 2, a line search is performed along the direction $s^{k,n}$ starting at $x^k$ with a step $\sigma_{k,n,l}$, evolving (or not) at each loop in $l$, until to $l$ reachs the value $j_{k,n}+1$, or until that the vectors $v^{k,n,l}$ and $\xi^{k,n,l}$ do not satisfy respectively the inequalities (\ref{vsstrict}) and (\ref{vxistrict}).

\item It is the evaluation of the pair $(y^{k,n,l},\xi^{k,n,l})$ and of $v^{k,n,l}$ at (2.g) which conducts to the step 3. The two possible ends at step 2 conduct either at a null step which increases the bundle, but keeps the iterate $x^k$ unchanged, or either to a serious step which increases also the bundle by returning at the step 1, and produces the new iterate $x^{k+1}$.
\end{itemize}

\begin{remark}\label{remapp}
We can make explicit the achievement of $s^{k,n,j}$ and $v^{k,n,j}$, obtained respectively at steps 1 and 2. According to the transportation formula, the convex hull generated by the vectors $w^i$, $i\in \widehat{I}_{k,n,j}$ is a polyhedral approximation of $T^{\hat{\epsilon}_{k,n,j}}(x^k)$, in the sense that for all $j$, the vector $$s^{k,n,j}:=\displaystyle\sum_{i\in \widehat{I}_{k,n,j}}\alpha_i^{k,n,j}w^i\in T^{\hat{\epsilon}_{k,n,j}}(\hat{x}^{k,n,j}),\mbox{ where }\hat{x}^{k,n,j}:=\displaystyle\sum_{i\in \widehat{I}_{k,n,j}}\alpha_i^{k,n,j}z^i\in\overline{B(x^{k},R\,2^{-j})},$$ with
\begin{eqnarray*}
    \hat{\epsilon}_{k,n,j} &=& \displaystyle\sum_{i\in \widehat{I}_{k,n,j}}\alpha_i^{k,n,j}\bigg\langle z^i-\sum_{m\in \widehat{I}_{k,n,j}}\alpha_m^{k,n,j}z^m,w^i-\sum_{m\in \widehat{I}_{k,n,j}}\alpha_m^{k,n,j}w^m\bigg\rangle\\
&\leq& 2\mu_{k,n,j}' R\,2^{-j}=\mu_{k,n,j}' R\,2^{-j+1},
\end{eqnarray*}
where $\mu'_{k,n,j}:=\max\{\|w^i\||i\in\widehat{I}_{k,n,j}\}$.

We can define in same way the vector $v^{k,n,l}$ belonging to the polyhedral approximation of $T^{\epsilon_{k,n,l}}(y^{k,n,l})$ by keeping only the $w^i$ whose arguments $z^i$ are in the neighbourhood $\overline{B(y^{k,n,l},R\,2^{-l})}$ of $y^{k,n,l}$. Let us justify that. At step (2.e), for all $l$, the vector $v^{k,n,l}$ is given by: $v^{k,n,l}=\displaystyle\sum_{i\in I_{k,n,l}}\lambda_iw^i$, with $\lambda\in\Delta_{I_{k,n,l}}$.
According to the transportation formula, one has:
$$v^{k,n,l}\in T^{\epsilon_{k,n,l}}(\hat{y}^{k,n,l}),\mbox{ where }\hat{y}^{k,n,l}:=\displaystyle\sum_{i\in I_{k,n,l}}\lambda_iz^i.$$

Then, Corollary \ref{corotf} gives that $$\hat{y}^{k,n,l}\in\overline{B(y^{k,n,l},R\,2^{-l})}\mbox{ and  }  
\epsilon_{k,n,l}\leq 2\mu_{k,n,l} R\,2^{-l}\leq\mu_{k,n,l} R\,2^{-l+1},$$
where $\mu_{k,n,l}:=\max\{\|w^i\||i\in I_{k,n,l}\}$.

Observe that in the case $\epsilon_k=0$, as we do not need a polyhedral approximation for $T^{\epsilon_k}(y^k)=T(y^k)$, we can take $v^k=\xi^k$. But in this case, the method HAEPPA brings back to the one developped in \cite{hppa}, and consequently, our algorithm of bundle type would be similar to the method in \cite{rbm}.
\end{remark}\bigbreak

\begin{remark}\label{rqpt4}
\begin{enumerate}
    \item Note that according to the definition of $\widehat{I}_{k,n,j}$, the pair $(x^k,u^k)\in\mathrm{Gr}(T)$ is always in the reduced bundle used to generate $s^{k,n,j}$ (because $\|x^k-x^k\|=0\leq R\,2^{-j}$), what guarantees that this last is nonempty, and then, that $\alpha^{k,n,j}$ is well defined. Likewise, the reduced bundle used to approach $T^{\epsilon_{k,n,l}}(y^{k,n,l})$, is nonempty because the set of index $I_{k,n,l}$ is also nonempty: it contains the element $i$ such that $z^i=x^k$; indeed, according to steps (2.b)-(2.c):
\begin{equation}\label{ineqykxkmaj}
 \|y^{k,n,l}-x^k\|= \sigma_{k,n,l}\|s^{k,n}\|\leq R\, 2^{-l}.   
\end{equation}
\bigbreak
    \item \label{iincichap} During the null step, one has one of the relations (\ref{vsstrict}) or/and (\ref{vxistrict}) with $l_{k,n}=j_{k,n}+1$. It follows that at the end of a null step, one has: $I_{k,n,l_{k,n}}\subset\widehat{I}_{k,n,j_{k,n}}$. 
    Indeed, for all $i\in I_{k,n,l_{k,n}}$:
\begin{eqnarray*}
\|z^i-x^k\| &\leq& \|z^i-y^{k,n}\|+\|y^{k,n}-x^k\|\\
&\leq& R\,2^{-l_{k,n}}+R\,2^{-l_{k,n}}\\
&=& 2R\,2^{-l_{k,n}}\\ &=& R\,2^{-j_{k,n}} \hspace{2cm}\mbox{because $l_{k,n}=j_{k,n}+1$.}
\end{eqnarray*}

Thus, $v^{k,n}\in\mathrm{conv}\{w^i\}_{i\in\widehat{I}_{k,n,j_{k,n}}}$.

    \bigbreak   

\item At the end of (3.a), the line search performed at step 2, having constructed a pair $(y^{k,n},\xi^{k,n})$ such that $\|y^{k,n}-x^k\|\leq R\,2^{-l}\leq R\,2^{-j-1}\leq R\,2^{-j}$ (because $l=j+1$), the pair $(y^{k,n},\xi^{k,n})$ is added to the sub-bundles $\{(z^i,w^i):i\in \widehat{I}_{k,n,j}\}$ by returning at step (1.b) (by recalling that $(y^{k,n},\xi^{k,n})$ is added to this reduced bundle if $\|y^{k,n}-x^k\|\leq R\,2^{-j}$).\bigbreak    

\item A serious step intervenes when $l_{k,n}\leq j_{k,n}+1$, with:
\begin{equation}\label{bigineq1}
    \langle v^{k,n},s^{k,n}\rangle \geq  \frac{1}{2}\|s^{k,n}\|^2,
\end{equation}
and 
\begin{equation*}
    \langle s^{k,n},\xi^{k,n}\rangle\geq \frac{1}{2}\|s^{k,n}\|^2. 
\end{equation*}
\end{enumerate}

\end{remark}\bigbreak

We show now that Algorithm \ref{bmhaeppa} verifies the scheme (\ref{vk2bm})-(\ref{ineqarror2bm}):
\begin{description}
    \item[a.] The relation (\ref{vk2bm}) is explained in Remark \ref{remapp}.
    \item[b.] The proximal relation (\ref{iterprox}) is due to the choice made for $y^k:=x^k-\sigma_ks^k$; indeed, by noting $e^k=c_kv^k-\sigma_ks^k$, one has:
\begin{eqnarray*}
    0 &=& \sigma_ks^{k}+(y^{k}-x^{k})\\
      &=& c_kv^{k}+(y^{k}-x^{k})-e^k.
\end{eqnarray*}
    \item[c.] The error's relation (\ref{ineqarror2bm}) is a consequence of the fact that each new iterate of the sequence $\{x^k\}$ comes from a serious step, and that at a serious step, it holds: 
\begin{equation}\label{geqmax}
    \langle v^{k},s^{k}\rangle \geq \frac{1}{2}\|s^{k}\|^2.
\end{equation}
Indeed, 
\begin{eqnarray*}
\|e^k\|^2 &=& \|c_kv^{k}+(y^{k}-x^{k})\|^2\\
          &=& \|c_kv^{k}-\sigma_ks^k\|^2\\
          &=& c_k^2\|v^k\|^2-2c_k\sigma_k\langle v^k,s^k\rangle+\sigma_k^2\|s^k\|^2\\
          &\leq& c_k^2\|v^k\|^2-c_k\sigma_k\|s^k\|^2+\sigma_k^2\|s^k\|^2\quad \mbox{because of (\ref{geqmax})}\\
          &\leq& c_k^2\|v^k\|^2+\sigma_k^2\|s^k\|^2\\
          &=& c_k^2\|v^k\|^2+\|y^k-x^k\|^2.
\end{eqnarray*}
    \item[d.] The relations (\ref{xikinct}) and (\ref{iteratexk2bm}) are explicit in the steps (2.c) and (3.b).
\end{description}
\bigbreak

\subsection{Convergence analysis}

We begin by giving a very important lemma, in the measure that it enables to prove the convergence of the algorithm when we replace $T^{\hat{\epsilon}_k}(x^k)$ and $T^{\epsilon_k}(y^k)$ by their respective polyhedral approximations.\bigbreak 

\begin{lem}\label{lemgamma}\rm \cite[Lemma IX.2.1.1]{MR1295240}\,\,\it
Let $\gamma>0$ fixed. Consider two infinite sequences $\{t^m\}$ and $\{\hat{t}^m\}$ satisfying for $m=1,2,...$:
\begin{equation}\label{condgamma}
    \langle t^i-t^{m+1},\hat{t}^m\rangle\geq\gamma\|\hat{t}^m\|^2,\quad \mathrm{for\,\, all}\,\,
 i=1,...,m.
\end{equation}
If $\{t^i\}$ is bounded, then $\hat{t}^m\rightarrow 0$ when $m\rightarrow\infty$.
\end{lem}\bigbreak

In the following, we will suppose that $T^{-1}(0)\neq\varnothing$.\bigbreak

Before proving the convergence of the sequence $\{x^k\}$, we need some preliminary technical results:\bigbreak

\begin{prop}\label{proptecres}
Let $x^k$ be the current iterate in Algorithm \ref{bmhaeppa}, and assume $x^k\notin T^{-1}(0)$. Then, after $x^{k+1}$ is generated  in (3.b), the following holds:
\begin{description}
\item[(i)] Let $H_{y^k,\xi^k}$ be the halfspace defined in (\ref{defhalfspace}), written with $(y,v):=(y^k,\xi^k)$. Then $x^k\notin H_{y^k,\xi^k}$ and $x^{k+1}=P_{H_{y^k,\xi^k}}(x^k)$.
\item[(ii)] For all $x^*\in T^{-1}(0)$, 
\begin{equation}\label{ineqxkdec2bm}
\|x^{k+1}-x^*\|^2\leq\|x^{k}-x^*\|^2-\|x^{k+1}-x^k\|^2.
\end{equation}
\item[(iii)] Finally, 
\begin{equation}\label{ineqxk1xktau}
    \|x^{k+1}-x^k\|>\tau\, R\,2^{-2(j_k+1)}/\|\xi^k\|.
\end{equation}
\end{description}
\end{prop}\bigbreak

\begin{proof}
(i) Recall that $x^k\in H_{y^k,\xi^k}$ if and only if $\langle x^k-y^k,\xi^k\rangle\leq 0$.
But, since $y^k=x^k-R\,2^{-l_k}s^k/\|s^k\|$, one has:
\begin{equation}\label{ineqxk1xkstrict}
\langle x^k-y^k,\xi^k\rangle=R\,2^{-l_k}\langle s^k,\xi^k\rangle/\|s^k\|\geq R\,2^{-l_k-1}\|s^k\|>\tau R\,2^{-l_k-j_k-1}>0,
\end{equation}
hence $x^k\notin H_{y^k,\xi^k}$. To see that $x^{k+1}$ is its projection, just recall the definition of $x^{k+1}$ in step (3.a) of Algorithm \ref{bmhaeppa}.\bigbreak

(ii) Combining the facts that $x^{k+1}$ is an orthogonal projection onto $H_{y^k,\xi^k}$ and $T^{-1}(0)\subset H_{y^k,\xi^k}$, (\ref{ineqxkdec2bm}) follows from the fact that orthogonal projections are firmly nonexpansive (see for example \cite{MR0388177}).\bigbreak

(iii) According to (\ref{ineqxk1xkstrict}), one has:
$$\|x^{k+1}-x^k\|=\langle x^k-y^k,\xi^k\rangle/\|\xi^k\|>\tau R\,2^{-l_k-j_k-1}/\|\xi^k\|\geq \tau R\,2^{-2j_k-2}/\|\xi^k\|.$$ 
\end{proof}\bigbreak

According to the item (ii) of Proposition \ref{proptecres}, the sequence $\{x^k\}$ generated by the  Algorithm \ref{bmhaeppa} is \textit{Fejér monotone relative to} $T^{-1}(0)$. This enables to show the boundedness of the variables in  the  Algorithm \ref{bmhaeppa}. The following result is essentially due to Opial \cite{MR0211301}:\bigbreak

\begin{prop}\label{propfejmonseq}
Consider a sequence $\{x^k\}$ satisfying for a nonempty set $\mathscr{S}$:
$$\|x^{k+1}-x\|\leq\|x^k-x\|,\quad for\, any\,\, x\in\mathscr{S}.$$
Such a sequence is said Fejér monotone relative to $\mathscr{S}$, and verifies:
\begin{description}
\item[(i)] $\{x^k\}$ is bounded.
\item[(ii)] If $\{x^k\}$ has a cluster point which is in $\mathscr{S}$, then the full sequence converges to a limit in $\mathscr{S}$.
\end{description} 
\end{prop}\bigbreak 

All the variables generated by Algorithm \ref{bmhaeppa}, that is $x^k$, $s^k$, $y^k$, $\xi^k$, $v^k$, $\{(y^{k,n},\xi^{k,n})\}$, $\{(z^p,w^p)\}$ and $\{v^{k,n}\}$ are bounded; in particular, one has the following result:\bigbreak

\begin{lem}\label{lemvarbor}
If $T^{-1}(0)\neq\varnothing$, the sequences $\{(y^{k,n},\xi^{k,n})\}$, $\{(z^p,w^p)\}$ and $\{v^{k,n}\}$  generated respectively by Algorithm \ref{bmhaeppa} at steps (2.h), (0.b)-(3.a) and (2.h), are bounded.
\end{lem}\bigbreak 

\begin{proof}
We show first that the sequence $\{x^k\}$ is bounded. If this sequence is finite, it is trivial; otherwise, $k\rightarrow\infty$. By Proposition \ref{proptecres}(ii), we deduce that the sequence $\{x^k\}$ is Fejér monotone relative to the nonempty set $T^{-1}(0)$. We can apply then Proposition \ref{propfejmonseq} to deduce that $\{x^k\}$ is bounded. Then, there exists a compact set $K_0$ such that the (bounded) sequence $\{x^k\}$ is included in $K_0$. Define $K_1:=K_0+\overline{B(0,R)}$; $K_1$ being also a compact set, $T$ locally bounded imply that $T(K_1)$ is bounded. According to the step (2.c), one has $y^{k,n,l}\in\overline{B(x^k,R2^{-l})}\subset K_1$: the sequence $\{y^{k,n,l}\}$ is then bounded. It follows from the steps (2.h) and (3.b), that the variables $y^{k,n}$ and $y^k$ are also bounded. As $\{z^p\}$ is extracted either of $\{x^k\}$, either of $\{y^k\}$, this sequence is also bounded. The sequences $\{\xi^{k,n}\}$ and $\{w^p\}$ being included in $T(K_1)$, they are bounded. As $\{v^{k,n}\}\subset\mathrm{conv}\{w^p\}$, $\{v^{k,n}\}$ is also bounded. 
\end{proof}\bigbreak
 
Before stating convergence results, we must show that our algorithm does not present infinite loops, and is consequently well defined.\bigbreak

The following result describe the different inner loops of the algorithm. We will retain that there is not infinite loop on those index when $x^k$ is not a zero of $T$.\bigbreak

\begin{lem}\label{longlem}
Let $x^k$ be the current iterate in Algorithm \ref{bmhaeppa}; suppose that $x^k\notin T^{-1}(0)$. Then:
\begin{description}
    \item[(i)] concerning the loop (1.e)$\leftrightarrow$(1.b), there exists a finite $j=j_{k,n}$ such that (1.f) is reached:
\begin{equation}\label{nors}
    \|s^{k,n}\|>\tau\, 2^{-j_{k,n}}.
\end{equation}
Furthermore, the loop (2.g)$\leftrightarrow$(2.b) is finite: (2.h) is reached with $l_{k,n}\leq j_{k,n}+1$.
    \item[(ii)] Concerning the loop (3.a)$\leftrightarrow$(1.b), there exists a finite $n=n_k$ such that (3.b) is reached.
\end{description}
\end{lem}\bigbreak

\begin{proof}	
Suppose that $0\notin T(x^k)$.
The transportation formula is used to show (i), that is, the loop in $j$ at step 1 is finite. 
Suppose for contradiction, that Algorithm \ref{bmhaeppa} loops forever in (1.e)$\leftrightarrow$(1.b). Then, $j\rightarrow\infty$ and an infinite sequence $\{s^{k,n,j}\}_{j\in\mathbb{N}}$ is generated, satisfying $\|s^{k,n,j}\|\leq\tau\, 2^{-j}$. Consequently, there exist subsequences $\{n_q\}$ and $\{j_q\}$ such that:
\begin{equation}\label{snqjq}
    \|s^{k,n_q,j_q}\|\leq\tau\, 2^{-j_q},
\end{equation}
with $\lim_{q\rightarrow\infty}j_q=\infty$. For such indices, define: $\widehat{I}_q:=\widehat{I}_{k,n_q,j_q}$. Because of step (1.b) one has: for all $i\in \widehat{I}_q$, $\|z^i-x^k\|\leq R\, 2^{-j_q}$. Consider the vector $\alpha^q:=\alpha^{k,n_q,j_q}$, which solves the problem of minimal norm at step  (1.c). Put:
\begin{equation*}
    (\hat{x}^q,\hat{s}^q):=\bigg(\sum_{i\in I_q}\alpha^q_iz^i,s^{k,n_q,j_q}\bigg).
\end{equation*}
Applying Corollary \ref{corotf} with $\rho=R\,2^{-j_q}$ and $\tilde{x}=x^k$, we obtain:
\begin{equation}\label{epdo}
    \hat{s}^q\in T^{\hat{\epsilon}_q}(\hat{x}^q),\quad \mathrm{with}\quad \hat{\epsilon}_q\leq 2R\,2^{-j_q}M,
\end{equation}
where $M:=\sup\Big\{\|u\||u\in T\big(\overline{B(x^k,R)}\big)\Big\}$. Furthermore, 
\begin{equation}\label{normeq}
    \|\hat{x}^q-x^k\|\leq R\,2^{-j_q}.
\end{equation}
Letting $q\rightarrow\infty$ in (\ref{snqjq}), (\ref{epdo}) and (\ref{normeq}), Proposition \ref{clteps} yields:
\begin{equation*}
    \lim_{q\rightarrow\infty}(\hat{\epsilon}_q,\hat{x}^q,\hat{s}^q\in T^{\hat{\epsilon}_q}(\hat{x}^q))=(0,x^k,0),
\end{equation*}
what implies that $0\in T^0(x^k)=T(x^k)$, and contradicts our starting hypothesis. Consequently, there exists a finite $j$ such that the loop (1.e)$\leftrightarrow$(1.b) ends, that is an indice $j_{k,n}$ such we have the inequality (\ref{nors}). We have then prove the first part of (i).\newline
The step 2 is always finite, because it suffies for example that $l=j_{k,n}+1$. Hence, the loop (2.g)$\leftrightarrow$(2.b) ends with the inequality $l=l_{k,n}\leq j_{k,n}+1$.\newline
Let us prove (ii) now. If an infinite loop occurs at (3.a)$\leftrightarrow$(1.b), then $n\rightarrow\infty$. Also, at step (2.h), an infinite sequence $\{(y^{k,n},v^{k,n})\}_{n\in\mathbb{N}}$ is generated. We have, for every $n$, the loop (1.e)$\leftrightarrow$(1.b) ends with an indice  $j$  such that:
\begin{equation}\label{doublecondition}
    s^{k,n,j}=s^{k,n,j_{k,n}}=s^{k,n}\quad\mathrm{et}\quad \|s^{k,n}\|>\tau\, 2^{-j_{k,n}}.
\end{equation}
Observe that in the loop (3.a)$\leftrightarrow$(1.b), the sequence $\{j_{k,n}\}_n$ is either incremented, or stay constant: in other words, this sequence is increasing (or even constant). If this sequence is not bounded, it would be divergent to infinity. Yet, 
it would be a contradiction with the assertion (i).
 Consequently, $j$ eventually reaches its final value, say $J$. Therefore, there exists $\bar{n}$ such that $j_{k,n}=j_{k,\bar{n}}=J$ for any $n\geq \bar{n}$. Consider now the sequence $\{(y^{k,n},v^{k,n})\}_{n\geq\bar{n}}$. At step (3.a), an infinity of null steps is made, and one has:
\begin{equation}\label{bigineq}
    \langle v^{k,n},s^{k,n}\rangle < \frac{1}{2}\|s^{k,n}\|^2
\end{equation}
or
\begin{equation}\label{strictvxi}
    \langle s^{k,n},\xi^{k,n}\rangle< \frac{1}{2}\|s^{k,n}\|^2, 
\end{equation}
with 
\begin{equation}\label{egalnorm}
    \|y^{k,n}-x^k\|\leq R\, 2^{-l_{k,n}}=R\, 2^{-j_{k,n}-1}=R\, 2^{-J-1},
\end{equation}
for all $n\geq\bar{n}$ (because if (2.a) leads to (3.a), then $l_{k,n}=j_{k,n}+1)$. This implies that $(y^{k,n}, \xi^{k,n})$ is incorporated to the sub-bundle associated to $\widehat{I}_{k,n,j_{k,n}}$ for any $n\geq\bar{n}$. In particular, if we choose an indice $\tilde{n}$ such that $\bar{n}<\tilde{n}<n$, then as $j_{k,\tilde{n}}=j_{k,n}=J$ in (\ref{egalnorm}), the pair $(y^{k,\tilde{n}},\xi^{k,\tilde{n}})$ is also incorporated in this sub-bundle defining $s^{k,n}$, which is the projection of 0 onto the convex hull of $\{w^i\}_{i\in \widehat{I}_{k,n,j_{k,n}}}$ (see steps (1.c)-(1.f)). A classical projection property gives:
\begin{equation}\label{ineqproj}
    \langle z-s^{k,n},-s^{k,n}\rangle\leq 0, \quad\forall z\in\mathrm{conv}\{w^i\}_{i\in \widehat{I}_{k,n,j_{k,n}}},
\end{equation}
with $z=\xi^{k,\tilde{n}}$, it can be written as: 
\begin{equation*}
    \langle \xi^{k,\tilde{n}}-s^{k,n},-s^{k,n}\rangle\leq 0, \quad\forall n>\tilde{n}>\bar{n},
\end{equation*}
hence:
\begin{equation}\label{ineqvs}
    \langle \xi^{k,\tilde{n}},s^{k,n}\rangle\geq \|s^{k,n}\|^2, \quad\forall n>\tilde{n}>\bar{n}.
\end{equation}
Suppose first that the inequality (\ref{strictvxi}) holds.
We have then:
\begin{equation*}
    -\langle \xi^{k,n},s^{k,n}\rangle>-\frac{1}{2}\|s^{k,n}\|^2. 
\end{equation*}
By suming this last inequality with the relation (\ref{ineqvs}), one has:
\begin{equation}\label{xitilde}
    \langle \xi^{k,\tilde{n}}-\xi^{k,n},s^{k,n}\rangle\geq (1-\frac{1}{2})\|s^{k,n}\|^2=\frac{1}{2}\|s^{k,n}\|^2, \quad\forall n>\tilde{n}>\bar{n}.
\end{equation}

Let $m\in\mathbb{N}$ such that $n=\bar{n}+m+1$. With this script, for any natural integers $\tilde{n}$ and $\bar{n}$ such that $n>\tilde{n}>\bar{n}$, it holds: $\tilde{n}\in\{\bar{n}+1,\bar{n}+2,...,\bar{n}+m\}=\{\bar{n}+i|i\in\{1,...,m\}\}$. Define then:
\begin{equation*}
    t^i:=\xi^{k,\bar{n}+i},\quad \hat{t}^i:=s^{k,\bar{n}+i+1}.
\end{equation*}
Using (\ref{xitilde}), it follows:
\begin{equation*}
    \langle t^i-t^{m+1},\hat{t}^m\rangle\geq \frac{1}{2}\|\hat{t}^m\|^2, \quad\forall i=1,...,m.
\end{equation*}
Then, the sequences $\{t^m\}$ and $\{\hat{t}^m\}$ satisfy the condition (\ref{condgamma}) of Lemma \ref{lemgamma}, with $\gamma=\displaystyle\frac{1}{2}>0$. Moreover, according to Lemma \ref{lemvarbor}, the sequence $\{t^i\}$ is bounded, and hence, the last part of this lemma is verified. That enables to conclude:
\begin{equation}\label{limt}
    \hat{t}^m\rightarrow 0\quad\mathrm{when}\quad m\rightarrow\infty.
\end{equation}
However, (\ref{doublecondition}) and the choice of $\bar{n}$ gives:
\begin{equation}\label{tmstrict}
    \|\tilde{t}^m\|=\|s^{k,\bar{n}+m+1}\|>\tau\,2^{-J}>0,
\end{equation}
what contradicts (\ref{limt}). Consequently, the loop (3.a)$\leftrightarrow$(1.b) must finish with a finite value of $n$, which proves  (ii) in the case where the inequality (\ref{strictvxi}) holds.\newline
Let us study now the case where the inequality (\ref{bigineq}) holds. Recall that according to the assertion \ref{iincichap} of Remark \ref{rqpt4}, $v^{k,\tilde{n}}\in\mathrm{conv}\{w^i\}_{i\in \widehat{I}_{k,\tilde{n},j_{k,\tilde{n}}}}$. Furthermore, as $j_{k,n}=j_{k,\tilde{n}}=J$, it follows $v^{k,\tilde{n}}\in\mathrm{conv}\{w^i\}_{i\in \widehat{I}_{k,n,j_{k,n}}}$. We can then rewrite the relation (\ref{ineqproj}) with $z=v^{k,\tilde{n}}$: 
\begin{equation}\label{ineqvs1}
    \langle v^{k,\tilde{n}},s^{k,n}\rangle\geq \|s^{k,n}\|^2, \quad\forall n>\tilde{n}>\bar{n}.
\end{equation}
Furthermore, the inequality (\ref{bigineq}) is equivalent to:
\begin{equation}\label{negequa}
    -\langle v^{k,n},s^{k,n}\rangle>-\frac{1}{2}\|s^{k,n}\|^2.
\end{equation}

Summing (\ref{ineqvs1}) and (\ref{negequa}), we have:
\begin{equation}\label{ineqb}
    \langle v^{k,\tilde{n}}-v^{k,n},s^{k,n}\rangle\geq \frac{1}{2}\|s^{k,n}\|^2, \quad\forall n>\tilde{n}>\bar{n}.
\end{equation}

With a same approach as previously, that is, by writing $n=\bar{n}+m+1$ ($m\in\mathbb{N}$), and by setting:
\begin{equation*}
    t^i:=v^{k,\bar{n}+i},\quad \tilde{t}^i:=s^{k,\bar{n}+i+1},
\end{equation*}
the relation (\ref{ineqb}) becomes:
\begin{equation}
    \langle t^i-t^{m+1},\tilde{t}^m\rangle\geq \frac{1}{2}\|\tilde{t}^m\|^2, \quad\forall i=1,...,m.
\end{equation}
Moreover, using the fact that $\{t^i\}$ is bounded (see Lemma \ref{lemvarbor}), we can apply again Lemma \ref{lemgamma} with $\gamma=\displaystyle\frac{1}{2}$. Then, (\ref{limt}) follows. We obtain a contradiction with the inequality (\ref{tmstrict}), and hence, there is a finite loop  between (3.a) and (1.b) in this case. The assertion (ii) is then proved.
\end{proof}\bigbreak

In the following result, we analyse the possibilities for an iteration of Algorithm \ref{bmhaeppa}:\bigbreak

\begin{prop}\label{2bmfinite}
Suppose that $T$ is a maximal monotone operator on $\mathbb{R}^N$, with $T^{-1}(0)\neq\varnothing$. Let $x^k$ be the current iterate in Algorithm \ref{bmhaeppa}. Then,
\begin{description}
    \item[(i)] if $x^k$ is a solution, either the oracle answers $u^k=0$ and the algorithm stops in (0.a), either the algorithm loops forever without updating $k$ (there is an infinity of null steps), after reaching the last serious step.
    \item[(ii)] Else, the algorithm reaches (3.b) after a finite number of inner iterations. Furthermore,
\begin{equation}\label{ineqskleq}
\|s^{k,n_k^*,j_k-1}\|\leq\tau 2^{-j_k+1},
\end{equation}
where $n_k^*$ is the smallest value of $n$ equating $j_{k,n}=j_k$, whenever $j_k>0$.
\end{description}  
\end{prop}\bigbreak

\begin{proof}
First, suppose that $x^k$ is a solution. In this case, if the oracle gives $u^k=0$, then the algorithm stops in (0.a). Otherwise, one has $u^k\neq 0$. Then, suppose by contradiction, that the step (3.b) is reached. We have:
\begin{eqnarray*}
    \langle x^k-y^k,\xi^k\rangle &=& \frac{\sigma_k}{2}\langle s^k,\xi^k\rangle\\
&=& R\, 2^{-l_k-1}\,\|s^k\|\\
&>& \tau\,R\, 2^{-(j_k+l_k+1)}>0, \quad\mathrm{with}\quad \xi^k\in T(y^k).
\end{eqnarray*}
But, as $0\in T(x^k)$, the last inequality contradicts the monotonicity of the operator $T$, that proves (i).\newline
Let us show now (ii). Suppose that $j_k>0$. If $x^k$ is not a solution, Lemma \ref{longlem} ensures that there is not infinite loop in $k$, and (3.b) is reached. 
To prove (\ref{ineqskleq}), define 
$$n_k^*:=\min\{n\leq n_k|j_{k,n}=j_k\}.$$
Then in (1.e)$\leftrightarrow$(1.b), the indices 
\begin{center}
$j:=j_{k,n_k^*}-1<j_k$ and $j+1=j_{k,n_k^*}=j_k$,
\end{center}
are such that (1.e) holds for index $j$ and (1.f) holds for index $j+1$:
\begin{center}
$\|s^{k,n_k^*,j_{k,n_k^*}-1}\|\leq\tau\, 2^{-(j_{k,n_k^*}-1)}$ and $\|s^{k,n_k^*,j_{k,n_k^*}}\|>\tau\, 2^{-j_{k,n_k^*}}$,
\end{center}
and the conclusion follows.
\end{proof}\bigbreak

This corollary follows immediately from the last result:\bigbreak

\begin{coro}\label{coro1}
The sequence of serious steps $\{x^k\}$ generated by Algorithm \ref{bmhaeppa} is either finite, with the last iterate, solution of the problem (\ref{inc1}); or either infinite, with none iterate solution (and in this second case, there exists always a finite number of null steps after each serious step).
\end{coro}\bigbreak

\begin{lem}\label{lemjk}
Suppose that Algorithm \ref{bmhaeppa} loops forever on $k$ (i.e. $k\rightarrow\infty$). Then $\lim_{k\rightarrow\infty}j_k=\infty$.
\end{lem}\bigbreak

\begin{proof}
According to the assertion (ii) of Proposition \ref{proptecres}, for any $x^*\in T^{-1}(0)$ one has:
$$\forall k\in\mathbb{N},\quad\|x^{k+1}-x^*\|^2\leq\|x^k-x^*\|^2.$$
We deduce then that the sequence $\{\|x^k-x^*\|\}$ is decreasing; furthermore, being bounded under by $0$, this sequence is convergent. It follows from the relation (\ref{ineqxkdec2bm}) that the sequence $\{\|x^{k+1}-x^k\|\}$ tends to $0$. 
Combining the fact that  $\{\|x^{k+1}-x^k\|\}$ tends to $0$ and $\{v^k\}$ is bounded in the relation (\ref{ineqxk1xktau}), we can conclude that  $\lim_{k\rightarrow\infty}j_k=\infty$. 
\end{proof}\bigbreak

We can sum up all the previous results with the following Theorem:\bigbreak

\begin{teo}
Consider the sequence $\{x^k\}$ generated by Algorithm \ref{bmhaeppa}. Then, the sequence is either finite with the last element in $T^{-1}(0)$, or it converges to a solution of (\ref{inc1}).
\end{teo}\bigbreak

\begin{proof}
We already dealt with the finite case in Proposition \ref{2bmfinite}.\\
If there are infinitely serious steps, as $\{x^k\}$ is Fejér monotone relative to $T^{-1}(0)$, we only need to show that one of its cluster point is a solution of (\ref{inc1}). The sequence $\{x^k\}$ being bounded, it admits a cluster point $x^*$, which in turn, is the limit of a subsequence $\{x^{k_q}\}$ of $\{x^k\}$. Because of Lemma \ref{lemjk}, we can deduce $j_{k_q}>0$ for $q$ large enough. Then, Proposition \ref{2bmfinite}(ii) applies: for $n_k^*$ defined therein, we have:
\begin{equation}\label{ineqskleqjkq}
\|s^{k_q,n_{k_q}^*,j_{k_q}-1}\|\leq\tau 2^{-j_{k_q}+1}.
\end{equation}
Consider the associated index set $I_q:=I_{k_q,n_{k_q}^*,j_{k_q}-1}$. By a same reasoning as in the proof of Lemma \ref{longlem}, mutatis mutandis, define 
\begin{eqnarray*}
\alpha^q &:=& \alpha^{k_q,n_{k_q}^*,j_{k_q}-1},\\
\hat{x}^q &:=& \sum_{i\in I_q}\alpha^q_iz^i,\\
\hat{s}^q &:=& s^{k_q,n_{k_q}^*,j_{k_q}-1}\,\,\,=\,\,\,\sum_{i\in I_q}\alpha^q_iw^i.
\end{eqnarray*}
We have that:
\begin{equation}\label{inqxqxkqleq}
\|\hat{x}^q-x^{k_q}\|\leq\sum_{i\in I_q}\alpha^q_i\|z^i-x^{k_q}\|\leq R\,2^{-(j_{kq}-1)}=R\,2^{-j_{kq}+1}.
\end{equation}
Let $M$ be a upper bound for $\|w^p\|$ (these variables are bounded by Lemma \ref{lemvarbor}). Then Corollary \ref{corotf} yields 
\begin{equation}\label{applilemtf}
\hat{s}^q\in T^{\hat{\epsilon}_q}(\hat{x}^q)\quad\mbox{with}\quad\hat{\epsilon}_q\leq 2 R\,2^{-j_{kq}+1}M.
\end{equation}
Using Lemma \ref{lemjk}, we have $\lim_{q\rightarrow\infty}j_{k_q}=\infty$. Hence, by (\ref{ineqskleqjkq}), (\ref{inqxqxkqleq}), (\ref{applilemtf}), we conclude that:
\begin{equation*}
(\hat{\epsilon}_q,\hat{x}_q,\hat{s}_q\in T^{\hat{\epsilon}_q}(\hat{x}^q))\longrightarrow(0,x^*,0),
\end{equation*} 
when $q\rightarrow\infty$. We can apply then Proposition \ref{clteps} to conclude that $0\in T(x^*)$.
\end{proof}\bigbreak

\section{Concluding Remarks}

A new bundle method for solving the inclusion problem for a maximal monotone operator $T$ on $\mathbb{R}^N$ was presented. This algorithm is the first presenting a double polyhedral approximation of the $\epsilon$-enlargement of a general maximal monotone operator. It could be a support for new bundle methods using also two polyhedral approximations obtained via the transportation formula; for example, we could create implementable versions of the methods in \cite{happa} and \cite{MR1871872} with a double polyhedral approximation. Also, as was pointed out by the referee, the issues concerning applicability of inexact proximal methods to  splitting/decomposition go way beyond forward-backward of \cite{mfb} (see \cite{MR2583900,MR2095352}). That furnishes other possibilities of adaptation with bundle methods.

\end{document}